\documentclass[12pt]{article}
\usepackage{graphicx}
\usepackage{amsmath,amsthm,amssymb,enumerate}
\usepackage{euscript,mathrsfs}
\usepackage[left=2cm,right=2cm,top=3.5cm,bottom=3.5cm]{geometry}
\usepackage{color}

\usepackage{soul}

\catcode`\@=11 \@addtoreset{equation}{section}

\catcode`\@=12

\allowdisplaybreaks

\newtheorem{Theorem}{Theorem}[section]
\newtheorem{Proposition}[Theorem]{Proposition}
\newtheorem{Lemma}[Theorem]{Lemma}
\newtheorem{Corollary}[Theorem]{Corollary}

\theoremstyle{definition}
\newtheorem{Definition}[Theorem]{Definition}

\newtheorem{Remark}[Theorem]{Remark}

\newcommand{\bTheorem}[1]{
\begin{Theorem} \label{T#1} }
\newcommand{\eT}{\end{Theorem}}

\newcommand{\bProposition}[1]{
\begin{Proposition} \label{P#1}}
\newcommand{\eP}{\end{Proposition}}

\newcommand{\bLemma}[1]{
\begin{Lemma} \label{L#1} }
\newcommand{\eL}{\end{Lemma}}

\newcommand{\bCorollary}[1]{
\begin{Corollary} \label{C#1} }
\newcommand{\eC}{\end{Corollary}}

\newcommand{\bRemark}[1]{
\begin{Remark} \label{R#1} }
\newcommand{\eR}{\end{Remark}}

\newcommand{\bDefinition}[1]{
\begin{Definition} \label{D#1} }
\newcommand{\eD}{\end{Definition}}

\newcommand{\Q}{\mathbb{T}^d}

\newcommand{\bfphi}{\boldsymbol{\varphi}}

\newcommand{\bFormula}[1]{
\begin{equation} \label{#1}}
\newcommand{\eF}{\end{equation}}

\newcommand{\DC}{C^\infty_c}

\newcommand{\vc}[1]{{\bf #1}}

\newcommand{\Div}{{\rm div}_x}
\newcommand{\Grad}{\nabla_x}

\newcommand{\dx}{\,{\rm d} {x}}

\newcommand{\dt}{\,{\rm d} t }

\newcommand{\vU}{\vc{U}}

\newcommand{\intQ}[1]{\int_{{\Q}} #1 \ \dx}

\newcommand{\vv}{\vc{v}}

\newcommand{\D}{{\rm d}}

\def\softd{{\leavevmode\setbox1=\hbox{d}%
          \hbox to 1.05\wd1{d\kern-0.4ex{\char039}\hss}}}
\definecolor{Cgrey}{rgb}{0.85,0.85,0.85}
\definecolor{Cblue}{rgb}{0.50,0.85,0.85}
\definecolor{Cred}{rgb}{1,0,0}
\definecolor{fancy}{rgb}{0.10,0.85,0.10}

\newcommand\Cbox[2]{%
    \newbox\contentbox%
    \newbox\bkgdbox%
    \setbox\contentbox\hbox to \hsize{%
        \vtop{
            \kern\columnsep
            \hbox to \hsize{%
                \kern\columnsep%
                \advance\hsize by -2\columnsep%
                \setlength{\textwidth}{\hsize}%
                \vbox{
                    \parskip=\baselineskip
                    \parindent=0bp
                    #2
                }%
                \kern\columnsep%
            }%
            \kern\columnsep%
        }%
    }%
    \setbox\bkgdbox\vbox{
        \color{#1}
        \hrule width  \wd\contentbox %
               height \ht\contentbox %
               depth  \dp\contentbox
        \color{black}
    }%
    \wd\bkgdbox=0bp%
    \vbox{\hbox to \hsize{\box\bkgdbox\box\contentbox}}%
    \vskip\baselineskip%
}

\date{}

\begin{document}


\title{On a class of generalized solutions to equations describing incompressible viscous fluids}

\author{Anna Abbatiello\thanks{The research of A.A. is supported by Einstein Foundation, Berlin.} \and Eduard Feireisl\thanks{The research of E.F. leading to these results has received funding from the
Czech Sciences Foundation (GA\v CR), Grant Agreement
18-12719S. The stay of E.F. at TU Berlin is supported by Einstein Foundation, Berlin.} }


\maketitle

\bigskip

\centerline{Faculty of Mathematics and Physics, Charles University}
\centerline{Sokolovsk\'a 83, CZ-186 75 Prague 8, Czech Republic}
\centerline{feireisl@math.cas.cz}
\centerline{and}

\centerline{Institute of Mathematics, Technische Universit\"{a}t Berlin,}
\centerline{Stra{\ss}e des 17. Juni 136, 10623 Berlin, Germany}
\centerline{anna.abbatiello@tu-berlin.de}

\begin{abstract}

We consider a class of viscous fluids with a general monotone dependence of the viscous stress on the symmetric velocity gradient.
We introduce the concept of \emph{dissipative solution} to the associated initial boundary value problem inspired by the measure--valued
solutions for the inviscid (Euler) system. We show the existence as well as the weak--strong uniqueness property in the class of dissipative solutions. Finally, the dissipative solution enjoying certain extra regularity coincides with a strong solution of the same problem.   

\end{abstract}


{\bf MSC:} 35 Q 35, 35 A 01, 35 A 02,

\section{Introduction}
\label{I}

The main goal of the present paper is to develop a mathematical theory of viscous fluids in the case 
of very low regularity of solutions of the underlying 
evolutionary equations, similar to the Euler system describing the inviscid fluids. In particular, we allow 
the viscous stress tensor to be merely bounded function of the velocity gradient, where the corresponding energy estimate provide only 
bounds of the total variation (in fact total ``deformation'') of the velocity. We refer to the recent survey by M\' alek, Blechta, and Rajagopal \cite{BleMalRaj} for the relevant physical 
background of the equations and systems studied below. 

The motion of a general viscous incompressible fluid is described in terms of its \emph{velocity} $\vv = \vv(t,x)$
satisfying the following system of equations:
\begin{equation} \label{i1}
\Div \vv = 0,
\end{equation}
\begin{equation} \label{i2}
\partial_t \vv + \Div (\vv \otimes \vv) + \Grad \Pi = \Div \mathbb{S}.
\end{equation}
Here $\Pi$ is the associated pressure and $\mathbb{S}$ denotes the viscous stress tensor related to the symmetric velocity gradient
\[
\mathbb{D} \vv \equiv \frac{1}{2} \left( \Grad \vv + \Grad^t \vv \right)  
\]
through a general ``implicit'' rheological law 
\begin{equation} \label{i3}
\mathbb{S}: \mathbb{D} \vv = F(\mathbb{D} \vv) + F^*(\mathbb{S})
\end{equation}
where 
\begin{equation} \label{i4}
F : R^{d \times d}_{\rm sym} \to [0, \infty), \ F(0) = 0,\ {\rm Dom}(F) = R^{d \times d}_{\rm sym},
\end{equation}
is a convex function and $F^*$ denotes its conjugate. Note that \eqref{i3} means 
\[
\mathbb{S} \in \partial F(\mathbb{D} \vv), \ \mbox{or, equivalently,}\ \mathbb{D} \vv \in \partial F^* (\mathbb{S}),
\] 
in particular, the mapping 
\[
\mathbb{D} \in R^{d \times d}_{\rm sym} \to \mathbb{S} \in \partial F(\mathbb{D}) \in R^{d \times d}_{\rm sym}
\]
is monotone. Here and hereafter, the symbol $\partial$ denotes the subdifferential of a convex function, while $ R^{d \times d}_{\rm sym}$ is the space of real symmetric tensors.

To avoid technicalities connected with the kinematic boundary, we restrict ourselves to the spatially periodic solutions defined on the flat torus 
\begin{equation} \label{i5}
\Q = \left( [-1,1]|_{\{ -1, 1 \}} \right)^d.
\end{equation}
The problem is formally closed by imposing the initial data 
\begin{equation} \label{i6}
\vv(0, \cdot) = \vv_0.
\end{equation}

\subsection{Weak and strong solutions}

Problem \eqref{i1}--\eqref{i6} is essentially well posed 
locally in time in the class of strong (classical) solutions for non--degenarate $F$, see e.g. Bothe and Pr\"{u}ss \cite{BotPru}. 
By non--degenerate we mean that the dependence of $\mathbb{S}$ on $\mathbb{D}$ is given as 
\begin{equation} \label{BoPr1}
\mathbb{S} = \mu ( |\mathbb{D} \vv|^2) \mathbb{D} \vv
\end{equation}
where $\mu$ is twice continuously differentiable function, 
\begin{equation} \label{BoPr2}
\mu(Z) > 0,\ 
\mu(Z) + 2 Z \mu'(Z) > 0 \ \mbox{for any}\ Z \in [0, \infty).
\end{equation}

An iconic example is the power law fluid, for which
\[
\mathbb{S} = \mu (|\mathbb{D} \vv|) \mathbb{D} \vv,\ 
\mu (|\mathbb{D}| ) = (\mu_1 + \mu_2 |\mathbb{D}|^2)^{\frac{p-2}{2}},\ p > 1,\  \mu_1\geq 0, \ \mu_2>0.  
\]
It is known that if $p \geq \frac{11}{5}$, and $d=3$, the system \eqref{i1}--\eqref{i6} possesses global in time strong solutions, 
see e.g. M\' alek, Ne\v cas, R\accent23u\v zi\v cka \cite{MNR}. The existence of weak solutions has been established in 
several cases, see e.g. the monograph by M\' alek et al. \cite{MNRR}. The best result  
dealing with the limit case $d = 3$, $p > \frac{6}{5}$ was obtained by Diening, R\accent23u\v zi\v cka, and Wolf \cite{DRW}.
More general rheological laws have been studied by Bul\' \i \v cek et al. \cite{BuGwMaSG}. 
The hypothesis $p > \frac{6}{5}$ represents a threshold in the 3D case. Indeed the energy balance (inequality) associated with 
\eqref{i1}, \eqref{i2} reads  
\begin{equation} \label{i7}
\frac{1}{2} \intQ{ |\vv|^2 (\tau, \cdot) } + 
\int_0^\tau \intQ{ \mathbb{S}(\mathbb{D} \vv) : \mathbb{D} \vv } \dt \leq 
\frac{1}{2} \intQ{ |\vv_0|^2 }.
\end{equation}
In view of the standard Sobolev embedding $W^{1,p} \hookrightarrow L^q$, $1 \leq q \leq \frac{3p}{3 - p}$, the solutions 
are \emph{a priori} bounded in the space $L^q$, $q > 2$ only if $p > \frac{6}{5}$. In the opposite case, the convective term 
$\vv \otimes \vv$ is only bounded in $L^\infty(0,T; L^1(\Q; R^{d \times d}_{\rm sym}))$, which makes the analysis rather delicate.
The theory of dissipative solutions proposed in this paper goes beyond this threshold considering the class of velocities 
for which $\mathbb{D} \vv$ is merely a measure. Possible concentrations that may be produced by the convective term are captured 
by a Reynolds viscous stress $\mathfrak{R}_v$ introduced below.     

\subsection{Dissipative solutions}
\label{D}

Inspired by \cite{BreFeiHof19}, we introduce the concept of \emph{dissipative solution} to problem \eqref{i1}--\eqref{i6}. 
The essential features of the approach are {\bf (i)} augmenting the family of unknowns by quantities that account for possible 
oscillations/concentrations, {\bf (ii)} considering the energy balance as an integral part of the definition of generalized solutions.

Accordingly, we shall deal with the following quantities:
\begin{itemize}
\item the velocity $\vv \in C_{{\rm weak}}([0,T]; L^2(\Q; R^d))$;
\item the viscous stress tensor $\mathbb{S} \in L^1(0,T; L^1(\Q; R^{d \times d}_{\rm sym}))$;
\item the Reynolds viscous stress $\mathfrak{R}_v$ that acounts for possible concentrations in the convective term, 
\[
\mathfrak{R}_v \in L^\infty(0,T; \mathcal{M}^+(\Q; R^{d \times d}_{{\rm sym}})).
\]

\end{itemize}

\begin{Remark} \label{RD1}

The symbol $\mathcal{M}^+(\Q; R^{d \times d}_{\rm sym})$ denotes the space of finite vector valued signed measures on $\Q$ ranging in the cone of positively definite symmetric matrices. Specifically, 
\[
\mathbb{M} \in \mathcal{M}^+(\Q; R^{d \times d}_{\rm sym}) \ \Leftrightarrow \ 
\int_{\Q} \varphi (\xi \otimes \xi): \D \mathbb{M} \geq 0 \ \mbox{for any}\ \xi \in R^d,\ 
\varphi \in \DC(\Q),\ \varphi \geq 0.
\]

\end{Remark}

The dissipative solutions, introduced in detail in Section \ref{DD} below, solve the following system of equations 
\begin{equation} \label{DD1}
\begin{split}
\Div \vv &= 0,\\
\partial_t \vv + \Div (\vv \otimes \vv) + \Grad \Pi &= \Div \mathbb{S} - \Div \mathfrak{R}_v,\ \vc{v}(0, \cdot) = \vv_0, \\
\frac{1}{2} \intQ{ | \vv|^2 (\tau, \cdot) } + \int_{\Q} \frac{1}{2} {\rm trace}[\mathfrak{R}_v](\tau)  &+ \int_0^\tau \intQ{ 
\Big[ F(\mathbb{D} \vv) + F^* (\mathbb{S}) \Big] }\dt   \leq \frac{1}{2} \intQ{ |\vv_0|^2 } 
\end{split}
\end{equation}
in the sense of distributions.

\begin{Remark} \label{RD1a}

Strictly speaking, the symmetric gradient $\mathbb{D} \vv$ will be merely a measure on $(0,T) \times \Q$. Accordingly, its composition 
$F(\mathbb{D} \vv)$ with a convex function $F$ must be interpreted in a generalized sense proposed by Demengel and Temam \cite{DemTemb}, \cite{DemTema}.

\end{Remark} 

Although formally underdetermined, the problem \eqref{DD1} enjoys the \emph{weak--strong uniqueness} property. The velocity $\vv$ associated to the dissipative solutions coincides with the velocity satisfying \eqref{i1}, \eqref{i2} in the classical sense as long as the latter exists. In such a case, the measure $\mathfrak{R}_v$ vanishes while  
\[
\mathbb{S} \in \partial F(\mathbb{D} \vv)
\ \mbox{for a.a.}\ (t,x) \in (0,T) \times \Q.
\]  

The rest of the paper is organized as follows:
\begin{itemize}
\item In Section \ref{DD}, we introduce the dissipative solutions and state our main results.
\item In Section \ref{E}, we show existence of a dissipative solution for a fairly general class of initial data.
\item In Section \ref{R}, we establish the relative energy inequality associated to system \eqref{DD1}, which represents a basic tool for proving stability results. In particular, we establish the weak--strong uniqueness property in Section \ref{WS}.
\item In Section \ref{A}, we show that dissipative solution that are sufficiently regular coincide with the standard strong solutions 
of the same problem.
\item Finally, we discuss possible extensions of the method in Section \ref{C}.
\end{itemize}

\section{Prelimaries, main results}
\label{DD}

We start by introducing the basic hypotheses imposed on the rheological relation \eqref{i3} between the viscous stress $\mathbb{S}$ and 
the symmetric part of the velocity gradient $\mathbb{D} \vv$. We recall Fenchel--Young inequality
\begin{equation} \label{DD2}
\mathbb{S} : \mathbb{D} \leq F(\mathbb{D}) + F^*(\mathbb{S}) \ \mbox{for any}\ \mathbb{D}, \mathbb{S} \in R^{d \times d}_{\rm sym}
\end{equation}
yielding 
\begin{equation} \label{DD3}
F^* (\mathbb{S}) = \sup_{\mathbb{D} \in R^{d \times d}_{\rm sym}} \left[ \mathbb{S} : \mathbb{D} - F(\mathbb{D}) \right], 
\ \mbox{in particular}\ F(0) = 0 \Rightarrow F^*(\mathbb{S}) \geq 0. 
\end{equation}
Moreover, in view of hypothesis \eqref{i4}, the domain of $F$ is the whole space $R^{d \times d}_{\rm sym}$ which implies that 
$\mathbb{F}^*$ is superlinear, specifically, 
\begin{equation} \label{DD4}
\liminf_{|\mathbb{S}| \to \infty} \frac{ F^*(\mathbb{S}) }{|\mathbb{S}|} = \infty.
\end{equation}

For technical reasons specified below, we will assume that the domain of $F^*$ contains a ball in $R^{d \times d}_{\rm sym}$: 
there is $r > 0$ such that 
\begin{equation} \label{DD5}
0 \leq F^*(\mathbb{S}) < \infty \ \mbox{for any}\ \mathbb{S} \in R^{d \times d}_{\rm sym},\ |\mathbb{S}| < r.
\end{equation}
This implies that $F$ grows at least linearly for large $\mathbb{D}$, 
\begin{equation} \label{DD6}
\liminf_{|\mathbb{D}| \to \infty} \frac{ F(\mathbb{D}) }{|\mathbb{D}|} > 0.
\end{equation}

In view of the constitutive relation \eqref{i3}, we may anticipate that generalized 
solutions of \eqref{i1}, \eqref{i2}, satisfying some form of the energy balance 
\eqref{i7}, will be fields $\vv \in C_{\rm weak}([0,T]; L^2(\Q; R^d))$ with bounded deformation 
\[
\mathbb{D} \vv \in \mathcal{M}((0,T) \times \Q; R^{d \times d}_{\rm sym}) 
\cap L^\infty(0,T; W^{-1,2}(\Q; R^{d \times d}_{\rm sym})), 
\]
where $\mathcal{M}((0,T) \times \Q; R^{d \times d}_{\rm sym})$ denotes the set of finite tensor--valued Radon measures.

Now, we may use the machinery developed by Demengel and Temam \cite{DemTema} to define 
\[
F(\mathbb{D} \vv) \in \mathcal{M}((0,T) \times \Q).
\]
To this end, a few technical hypotheses are needed. Let
\[
F_\infty (\mathbb{D}) \equiv \lim_{s \to \infty} \frac{F(s \mathbb{D})}{s} \in [0, \infty]
\] 
be the asymptotic function of $F$. Following Demengel and Temam \cite{DemTema}, we further suppose that
there is $r > 0$ such that  
\begin{equation} \label{DD7}
{\rm Dom}(F^*) \subset \Lambda^0 + B(r,0), \ \mbox{where}\ 
\Lambda \equiv {\rm Dom}(F_\infty).
\end{equation}
Here $B(r,0)$ is the ball of radius $r$ centered at zero, while $\Lambda^0$ is the polar set of $\Lambda$ and ${\rm Dom}(\cdot)$ denotes the domain of  a function. In view of \cite[Proposition 1.2]{DemTema}, hypothesis \eqref{DD7} is equivalent to 
\[
F_{\infty} (\mathbb{D}) \leq r |\mathbb{D}| \ \mbox{for all}\ \mathbb{D} \in {\rm Dom}(F_\infty).
\]

Under these circumstances, one can define 
\[
F(\mathbb{D} \vv) \in \mathcal{M}((0,T) \times \Q),
\]
see \cite[Section 2]{DemTema}.

\begin{Remark} \label{RRDD1}

The hypotheses \eqref{DD5}, \eqref{DD7} may seem rather awkward at first glance. We claim they are automatically satisfied if 
${\rm Dom}(F^*) = R^{d \times d}_{\rm sym}$, in which case $F$ is superlinear in $\mathbb{D}$. We refer to 
\cite{DemTemb}, \cite{DemTema} for other interesting examples.

\end{Remark}

\subsection{Dissipative solutions}

As already pointed out, the dissipative solutions are defined in terms of the velocity field $\vv$, 
the viscous stress $\mathbb{S}$, and
the Reynolds viscous stress tensor $\mathfrak{R}_v$.

\begin{Definition} \label{DDD1}

We say that $\vv \in C_{\rm weak}([0,T]; L^2(\Q; R^d))$ is a \emph{dissipative solution} of the problem 
\eqref{i1}--\eqref{i6} if 
\begin{equation} \label{D33}
\mathbb{D} \vv \in \mathcal{M}((0,T) \times \Q; R^{d \times d}_{\rm sym}),
\end{equation}
and
there exist 
\[
\mathbb{S} \in L^1((0,T) \times \Q; R^{d \times d}_{\rm sym}), \ 
\mathfrak{R}_v \in L^\infty(0,T; \mathcal{M}^+(\Q; R^{d \times d}_{\rm sym}))
\]
such that the following holds: 

\begin{itemize}

\item {\bf Incompressibility.}
\begin{equation} \label{D3}
\int_0^T \intQ{ \vv \cdot \Grad \varphi } \dt = 0 
\end{equation}
for any $\varphi \in C^1([0,T] \times \Q)$.
\item {\bf Momentum equation.}
\begin{equation} \label{D4}
\begin{split}
&\intQ{ \vv \cdot \bfphi (\tau, \cdot) } - \intQ{ \vv_0 \cdot \bfphi(0, \cdot) } \\
&=\int_0^\tau\!\!\! \intQ{ \left[ \vv \cdot \partial_t \bfphi + (\vv \otimes \vv)\!:\! \Grad \bfphi \right] \!} \dt 
- \int_0^\tau\!\!\! \intQ{ \mathbb{S} \!:\! \Grad \bfphi \!} \dt + \int_0^\tau\!\! \int_{\Q} \Grad \bfphi \!:\! \D \mathfrak{R}_v \dt
\end{split}
\end{equation}
for any $0 \leq \tau \leq T$, and $\bfphi \in C^1([0,T] \times \Q; R^d)$, $\Div \bfphi = 0$.
\item{\bf Energy inequality.}
\begin{equation} \label{D5}
\intQ{ \frac{1}{2} |\vv|^2 (\tau, \cdot) } + \int_{\Q} \D \frac{1}{2}{\rm trace}[\mathfrak{R}_v (\tau)]
+ \int_0^\tau \!\!\intQ{\! \Big[ F(\mathbb{D} \vv) + F^*(\mathbb{S}) \Big] \!\!} \dt  
\leq \intQ{ \frac{1}{2} |\vv_0|^2 }
\end{equation}
for a.a. $\tau \in (0,T)$.

\end{itemize}

\end{Definition}

\begin{Remark} \label{RD7}

The integral $\int_0^\tau \intQ{ F(\mathbb{D} \vv) }$ should be understood as 
\[
\int_0^\tau \intQ{ F(\mathbb{D} \vv) } = \int_{(0,\tau) \times \Q} \D F(\mathbb{D} \vv)
\]
in the sense of Demengel and Temam \cite{DemTema}. Note that this extended definition is only necessary if 
${\rm Dom}(F^*)$ is not the whole space $R^{d \times d}_{\rm sym}$, because in such case one can prove that $F$ is superlinear and $\mathbb{D} \vv\in L^1((0,T);L^1(\Q, R^{d\times d}))$. 

\end{Remark}

\subsection{Main results}

As the first result, we state the existence of dissipative solutions proved in Section \ref{E}. 

\begin{Theorem} \label{TM1}

Let $\vv_0 \in L^2(\Q; R^d)$, $\Div \vv_0 = 0$, $d=2,3$, and $T > 0$ be given. Suppose that $F$ satisfies \eqref{i4}, together with 
\eqref{DD5}, \eqref{DD7}.

Then the problem \eqref{i1}--\eqref{i6} admits a dissipative solution $\vv \in C_{\rm weak}([0,T]; 
L^2(\Q; R^d))$ in the sense of Definition \ref{DDD1}.

\end{Theorem}

\begin{Remark} \label{RD2}

To simplify presentation, we deliberately omitted the issue of boundary conditions. However, the proof can be easily adapted to include some of the standard boundary conditions, in particular the complete slip or the no--slip boundary conditions. Note that in the latter case, the viscous stress $\mathbb{S}$ must be non--degenerate so that the trace of $\vv$ is well--defined. We refer to Blechta, M\' alek and Rajagopal \cite{BleMalRaj} for a detailed discussion of various boundary conditions.

\end{Remark}

The second result concerns the weak--strong uniqueness property. In view of the theory developed by Bothe and Pr\"{u}ss \cite{BotPru}, we state it in the $L^p-$framework.  

\begin{Theorem} \label{TM2}

Let $F$ satisfy the hypotheses \eqref{i4}, \eqref{DD5}, \eqref{DD7}.
Let $p > d + 2$, $d=2,3$. Suppose that the system \eqref{i1}--\eqref{i6} admits a strong solution $\widehat{\vv}$ in the class 
\begin{equation} \label{reg}
\widehat{\vv} \in C([0,T]; W^{2 - \frac{2}{p},p}(\Q; R^d)),\ 
\partial_t \widehat{\vv} \in L^p(0,T; L^p(\Q; R^d)),\ 
\widehat{\vv} \in L^p(0,T; W^{2,p}(\Q; R^d)) 
\end{equation}
defined on a time interval $[0,T]$, with the viscous stress 
\begin{equation}\label{new}
\widehat{\mathbb{S}}\in L^p(0,T; W^{1,p}(\Q, R^{d\times d}_{\rm sym})) \cap C([0,T]\times\Q, R^{d\times d}_{\rm sym}),
\end{equation}
and with the initial datum
\[
\widehat{\vv}(0, \cdot) = \vv_0 \in W^{2 - \frac{2}{p},p}(\Q; R^d).
\] 

Let $\vv$, together with $\mathbb{S}$ and $\mathfrak{R}_v$, be a dissipative solution of the same problem starting from the same initial data 
\[
\vv(0, \cdot) = \vv_0.
\]

Then 
\[
\mathfrak{R}_v = 0,\ \mathbb{S}(t,x) \in \partial F (\mathbb{D} \widehat{\vv} (t,x))  \ \mbox{for a.a.}\ (t,x) \in (0,T) \times \Q,
\]
and
\[
\vv (t,x) = \widehat{\vv}(t,x) \ \mbox{for a.a.}\ (t,x) \in (0,T) \times \Q.
\]
\end{Theorem}

\begin{Remark} \label{RD3}

Note that for $p > d  + 2$, we have the embedding relation $W^{2 - \frac{2}{p},p}(\Q) \hookrightarrow 
C^1(\Q)$. 

\end{Remark}

The proof of Theorem \ref{TM2} will be given in Sections \ref{R}, \ref{WS}. Note that \emph{existence} of 
local--in--time strong solutions in the class specified 
in Theorem \ref{TM2} was proved by Bothe and Pr\"{u}ss \cite{BotPru} for smooth non--degenerate viscous stress 
tensor satisfying \eqref{BoPr1}, \eqref{BoPr2}. In this case case \eqref{new} follows from \eqref{reg}.

Finally, we claim that dissipative solutions enjoying certain regularity are in fact strong solutions.

\begin{Theorem} \label{TM3}

Let $F$ satisfy the hypotheses \eqref{i4}, \eqref{DD5}, \eqref{DD7}.
Let $p > d + 2$, $d=2,3$. Let 
\[
\vv_0 \in W^{2 - \frac{2}{p},p}(\Q; R^d), \ \Div \vv_0 = 0
\]
be given. Suppose that $\vv$ is a dissipative solution on the system \eqref{i1}--\eqref{i6} on the time interval $[0,T]$
belonging to the regularity class  
\begin{equation} \label{reg1}
{\vv} \in C([0,T]; W^{2 - \frac{2}{p},p}(\Q; R^d)),\ 
\partial_t {\vv} \in L^p(0,T; L^p(\Q; R^d)),\ 
{\vv} \in L^p(0,T; W^{2,p}(\Q; R^d)).
\end{equation}

Then $\mathfrak{R}_v = 0$ and 
\[
\mathbb{S} : \mathbb{D} \vv = F(\mathbb{D} \vv) + F^*(\mathbb{S}),
\]
meaning $\vv$ is a strong solution of the system \eqref{i1}--\eqref{i6}.

\end{Theorem}

The proof of Theorem \ref{TM3} is given in Section \ref{A}.
\section{Existence}
\label{E}

For any $n\in \mathbb{N}$ we look for approximations $\vv^n$ and $\mathbb{S}^n$ that are solutions of the following system 
\begin{subequations}\begin{align} \label{galerkin}
\intQ{\!\!\partial_t \vv^n\! \!\cdot \bfphi_i\! } \!
 +\! \intQ{\!\! \mathbb{S}^n\! :\! \mathbb{D} \bfphi_i\!} 
&=\!\intQ{\!\! \vv^n\! \otimes \!\vv^n\!: \!\Grad \bfphi_i\!\!}, \ i=1,\dots,n,  \mbox{ a.e. in } (0,T), \\ 
\vv^n(0, \cdot)&=P^n\vv_0,
\end{align}\end{subequations}
such that they satisfy
\begin{equation}\label{subdif}
\mathbb{S}^n\in \partial F(\mathbb{D}\vv^n),
\end{equation}
and the following energy balance holds
\begin{equation}\label{eneq}
\frac{1}{2} \intQ{ |\vv^n(\tau) |^2 } + \int_0^\tau\!\! \intQ{ F^*(\mathbb{S}^n)+ F(\mathbb{D} \vv^n)\!} \dt =  \frac{1}{2}
\intQ{ |P^n \vv_0|^2 }\  \mbox{ for any}\ \tau\in (0,T).
\end{equation}
Here, the sequence $\{\bfphi_i\}_{i \in\mathbb{N}}$ denote the eigenfunctions of the Stokes operator, while  $P^n$ is the projection onto the finite dimensional space $[\bfphi_1,\dots,\bfphi_n]$ for any $n\in \mathbb{N}$.
The existence of such  $\vv^n$ and $\mathbb{S}^n$ can be achieved regularizing $\mathbb{S}^n$ by the introduction of a convolution kernel depending on $\mathbb{D}\vv^n$. With this approximation scheme the Galerkin method can be applied, then $\vv^n$ and $\mathbb{S}^n$ are obtained as limit in the finite dimensional spaces. A detailed proof was carried over in \cite[Section 3.1]{BuGwMaSG}. 

The energy inequality \eqref{eneq} together with \eqref{subdif} imply the following
\begin{align}\label{star}
&\vv^n\rightharpoonup^*\vv \mbox{ in } L^\infty(0, T; L^2(\Q; R^d)), \mbox{ as } n \to +\infty, \\
&\mathbb{S}^n\rightharpoonup  \mathbb{S}  \mbox{ in }L^1(0, T; L^1(\Q; R^{d\times d})), \mbox{ as } n \to +\infty, \label{convS}\\
&\mathbb{D}\vv^n\rightharpoonup^*\mathbb{D}\vv \mbox{ in }  \mathcal{M}((0,T) \times \Q; R^{d \times d}_{\rm sym}),  \mbox{ as } n \to +\infty,\label{convD}
\end{align}
at least for suitable subsequences, which we do not relabel. Note that \eqref{convS} follows from the superlinearity of $F^*$ and de la Vall\'{e}-Poussin criterion, while for \eqref{convD}, we use the assumption \eqref{DD6} and the embedding $L^1((0,T) \times \Q; R^{d \times d}_{\rm sym})\hookrightarrow \mathcal{M}((0,T) \times \Q; R^{d \times d}_{\rm sym})$.

As a consequence of \eqref{star}, we have 
$$ 0=\int_0^T\!\! \intQ{ \!\vv^n\! \cdot\! \Grad \varphi \!} \dt \to \int_0^T \!\!\intQ{\! \vv \cdot \Grad \varphi \!} \dt \mbox{ as } n \to +\infty,$$
for any $\varphi \in C^1([0,T] \times \Q)$, thus the incompressibility condition \eqref{D3} is satisfied.
 
Now, consider $\vv^n\otimes\vv^n-\vv\otimes\vv$. In view of the energy inequality \eqref{eneq}, the sequence is uniformly bounded in the space $L^\infty(0,T; L^1(\Q; R^{d\times d}))$, which is embedded into $L^\infty(0,T; \mathcal{M}(\Q; R^{d\times d}))$. Therefore there exists $\mathfrak{R}_v \in L^\infty(0,T; \mathcal{M}(\Q; R^{d \times d}_{\rm sym}))$ such that
\begin{equation}\label{abb}
(\vv^n\otimes\vv^n-\vv\otimes\vv) \rightharpoonup^*\mathfrak{R}_v \mbox{ in } L^\infty(0,T; \mathcal{M}(\Q; R^{d \times d}_{\rm sym})) \mbox{ as } n \to +\infty. 
\end{equation}
Now, through the same computations as in \cite[Section 3.2]{MarEd} it can be showed that $\mathfrak{R}_v$ is positive semidefinite, thus 
$$ \mathfrak{R}_v \in L^\infty(0,T; \mathcal{M}^+(\Q; R^{d \times d}_{\rm sym})).$$
Finally, from \eqref{galerkin} fixing $\bfphi\in [\bfphi_1,\dots, \bfphi_k]$ with $k\leq n$, multiplying by $\eta\in C^1(0,T)$, adding and subtracting $\int_0^t\intQ{\vv\otimes\vv:\nabla_x\bfphi\!\!}\ \eta(\tau)\dt$, taking the limit as $n\to +\infty$ and collecting the established convergences, using also a standard density argument, one gets \eqref{D4}.

In order to prove \eqref{D5} we rewrite \eqref{eneq} in the weak ``differential'' form: 
\begin{equation} \label{abb3}
- \int_0^T \partial_t \psi \intQ{ \frac{1}{2} |\vv^n|^2 } \dt + \int_0^T \psi \intQ{ \Big[ F(\mathbb{D} \vv^n) + 
F^* (\mathbb{S}^n ) \Big] } \dt = \intQ{ \frac{1}{2} |P^n \vv_0 |^2 } 
\end{equation}
for any $\psi \in \DC[0, T)$, $\psi(0) = 1$, $\partial_t \psi \leq 0$. Next, we write 
\[
\frac{1}{2} |\vv^n|^2 = \frac{1}{2} {\rm trace} \left[ \vv^n \otimes \vv^n - \vv \otimes \vv \right] + \frac{1}{2} |\vv|^2.
\]

Employing \cite[Lemma 3.1]{DemTema} we obtain
 \begin{equation} \label{abb4}
\int_0^T \int_{\Q} \widetilde{\psi} \D F(\mathbb{D} \vv) \equiv \int_0^T  \intQ{ \widetilde{\psi}  F(\mathbb{D} \vv)\!} \dt\leq \liminf_{n\to +\infty} \! \int_0^T\!\! \intQ{ \widetilde{\psi} F(\mathbb{D} \vv^n)\!} \dt
\end{equation}
for any $\widetilde{\psi} \in \DC(0,T)$, $0 \leq \widetilde \psi \leq \psi$,
where the measure $F(\mathbb{D} \vv)$ is understood in the sense of Demengel and Temam \cite{DemTema}. Moreover, 
since $F^*$ is a lower semicontinuous convex function, then
 \begin{equation} \label{abb5}
 \int_0^T \psi \intQ{\!\! F^*(\mathbb{S})\!} \dt\leq \liminf_{n\to +\infty} \! \int_0^T \psi \intQ{\!\!F^*(\mathbb{S}^n)\!} \dt.
\end{equation}
Consequently, using \eqref{abb}, \eqref{abb4}, and \eqref{abb5} we may perform the limit 
$n\to +\infty$ in \eqref{abb3} to obtain \eqref{D5}. This completes the proof of Theorem \ref{TM1}.

\section{Relative energy and weak--strong uniqueness}
\label{R}

We show the relative energy inequality that is a crucial tool for proving the weak strong uniqueness property
claimed in Theorem \ref{TM2}.

\subsection{Relative energy inequality}
  
Let $\vU$ be a continuously differentiable vector field, $\Div \vU = 0$.
It follows from \eqref{D4}, \eqref{D5} by a straightforward manipulation that 
\[
\begin{split}
\left[ \intQ{ \frac{1}{2} |\vv - \vU|^2 } \right]_{t = 0}^{t = \tau} &+ 
\int_{\Q} \D \frac{1}{2}{\rm trace}[\mathfrak{R}_v (\tau)]
+ \int_0^\tau \intQ{ \left[ F(\mathbb{D} \vv) + F^* (\mathbb{S}) \right] } \dt
\\
&\leq - \int_0^\tau \intQ{ \left[ \vv \cdot \partial_t \vU + (\vv \otimes \vv): \Grad \vU \right] } \dt 
+ \int_0^\tau \intQ{ \mathbb{S} : \Grad \vU } \dt \\ &- \int_0^\tau \int_{\Q} \Grad \vU : \D \mathfrak{R}_v \dt
+ \int_0^\tau \intQ{ \partial_t \vU \cdot \vU } \dt.
\end{split}
\] 
Indeed this follows by considering $\bfphi = \vU$ as test function in the momentum balance \eqref{D4} and summing up the resulting 
expression with the energy inequality \eqref{D5}.

Thus, regrouping some terms, we get 
\begin{equation} \label{R1}
\begin{split}
\left[ \intQ{ \frac{1}{2} |\vv - \vU|^2 } \right]_{t = 0}^{t = \tau} &+ 
\int_{\Q} \D \frac{1}{2}{\rm trace}[\mathfrak{R}_v (\tau)]
+ \int_0^\tau \intQ{ \left[ F(\mathbb{D} \vv) + F^*(\mathbb{S}) - \mathbb{S} : \mathbb{D} \vU \right] } \dt
\\
&\leq \int_0^\tau \intQ{ \left[ (\vU - \vv) \cdot \partial_t \vU - (\vv \otimes \vv): \Grad \vU \right] } \dt 
\\ &- \int_0^\tau \int_{\Q} \Grad \vU : \D \mathfrak{R}_v \dt.
\end{split}
\end{equation}
Finally, 
\[
\begin{split}
\int_0^\tau &\intQ{ (\vv \otimes \vv): \Grad \vU } \dt = 
\int_0^\tau \intQ{ \vv \otimes (\vv - \vU): \Grad \vU } \dt + \int_0^\tau \intQ{ (\vv \otimes \vU): \Grad \vU } \dt\\
&=  \int_0^\tau \intQ{ ( \vv - \vU)  \otimes (\vv - \vU): \Grad \vU } \dt 
+ \int_0^\tau \intQ{ \vU \otimes (\vv - \vU): \Grad \vU } \dt,
\end{split}
\]
where we have used that
\[
\int_0^\tau \intQ{ (\vv \otimes \vU): \Grad \vU } \dt = 0
\]
thanks to the incompressibility constraint \eqref{D3}.

Accordingly, relation \eqref{R1} reduces to 
\begin{equation} \label{R2}
\begin{split}
\left[ \intQ{ \frac{1}{2} |\vv - \vU|^2 } \right]_{t = 0}^{t = \tau} &+ 
\int_{\Q} \D \frac{1}{2}{\rm trace}[\mathfrak{R}_v (\tau)]
+ \int_0^\tau \intQ{ \Big[ F(\mathbb{D} \vv) + F^* (\mathbb{S}) - \mathbb{S}: \mathbb{D} \vU \Big] } \dt
\\
&\leq \int_0^\tau \intQ{ \left[ (\vU - \vv) \cdot \left( \partial_t \vU + \vU \cdot \Grad \vU \right) \right] } \dt 
\\ &- \int_0^\tau \int_{\Q} \Grad \vU : \D \mathfrak{R}_v \dt- 
\int_0^\tau \intQ{ ( \vv - \vU)  \otimes (\vv - \vU): \Grad \vU } \dt 
\end{split}
\end{equation}
Relation \eqref{R2} is called \emph{relative energy inequality}. Its validity can be easily extended by density argument to any vector field 
$\vU$ belonging to the regularity class \eqref{reg} from Theorem \ref{TM2}.

\section{Weak strong uniqueness property}
\label{WS}

The weak strong uniqueness follows immediately from the relative energy inequality \eqref{R2} applied to $\vU 
= \widehat{\vv}$ - the strong solution emanating from the initial data $\vv_0$. Indeed we have 
\[
\int_0^\tau \intQ{ (\widehat{\vv} - \vv) \cdot (\partial_t \widehat{\vv} + \widehat{\vv} \cdot \Grad \widehat{\vv}) } \dt = 
\int_0^\tau \intQ{ (\widehat{\vv} - \vv) \cdot \Div \widehat{\mathbb{S}} } \dt,  
\]
where, in accordance with \eqref{i3},  
\[
\widehat{\mathbb{S}} : \mathbb{D} \widehat{\vv} = F(\mathbb{D} \widehat{\vv}) + F^* (\widehat{\mathbb{S}}).
\]

Consequently, the relative energy inequality \eqref{R2} yields 
\begin{equation} \label{R3}
\begin{split}
\intQ{ \frac{1}{2} |\vv - \widehat{\vv} |^2 (\tau, \cdot) }  &+ 
\int_{\Q} \D \frac{1}{2}{\rm trace}[\mathfrak{R}_v (\tau)]
+ \int_0^\tau \intQ{ \Big[ F(\mathbb{D}\vv) + F^*(\mathbb{S}) - \mathbb{S} : \mathbb{D} \widehat{\vv}  \Big] } \dt\\
&+ \int_0^\tau \intQ{ (\vv - \widehat{\vv}) \cdot \Div \widehat{\mathbb{S}} } \dt
\\ &\leq - \int_0^\tau \int_{\Q} \Grad \widehat{\vv} : \D \mathfrak{R}_v \dt- 
\int_0^\tau \intQ{ ( \vv - \widehat{\vv})  \otimes (\vv - \widehat{\vv}): \Grad \widehat{\vv} } \dt 
\end{split}
\end{equation}
for a.a. $\tau \in (0,T)$.

Strictly speaking, $F(\mathbb{D}\vv)$ is a measure defined on an open set $(0, \tau)\times \Q$. Accordingly we should interpret
$$ \int_0^\tau\!\int_{\Q} F(\mathbb{D}\vv)dx\dt\equiv \int_0^\tau\!\int_{\Q} dF(\mathbb{D}\vv)\dt =\lim_{n\to+\infty}  \int_0^\tau\! \psi_n\int_{\Q}  dF(\mathbb{D}\vv)\dt$$
where $\psi_n\in C_c^\infty(0,\tau]$, $\psi_n \nearrow 1_{(0, \tau)}$, $0\leq \psi_n\leq 1$.
\\
Similarly,
we have 
\begin{align*}
&\int_0^\tau  \intQ{ (\vv - \widehat{\vv}) \cdot \Div \widehat{\mathbb{S}} } \dt= \lim_{n\to+\infty} 
\int_0^\tau \psi_n (t) \intQ{ (\vv - \widehat{\vv}) \cdot \Div \widehat{\mathbb{S}} } \dt\\
&= - \lim_{n\to+\infty} \int_0^\tau \psi_n (t) \intQ{ (\mathbb{D} \vv - \mathbb{D} \widehat{\vv}) : \widehat{\mathbb{S}} } \dt=  \int_0^\tau \intQ{\widehat{\mathbb{S}}:\mathbb{D} \widehat{\vv}}-\int_0^\tau \int_{\Q} \widehat{\mathbb{S}}:d\mathbb{D} \vv.  
\end{align*}

Thus the integrals on the left--hand side of \eqref{R3} can be handled by means of Fenchel--Young inequality:
\[
\begin{split}
\int_0^\tau &\intQ{ \Big[ F(\mathbb{D}\vv) + F^*(\mathbb{S}) - \mathbb{S} : \mathbb{D} \widehat{\vv}  \Big] } \dt
+ \int_0^\tau \intQ{ (\vv - \widehat{\vv}) \cdot \Div \widehat{\mathbb{S}} } \dt\\
&= \int_0^\tau \intQ{ \Big[ F(\mathbb{D}\vv) - \widehat{\mathbb{S}}:(\mathbb{D}\vv - \mathbb{D} \widehat{\vv}) - 
F(\mathbb{D} \widehat{\vv}) \Big] } \dt\\ 
&+ \int_0^\tau \intQ{ \Big[ F(\mathbb{D} \widehat{\vv}) + F^*(\mathbb{S}) - \mathbb{S} : \mathbb{D} \widehat{\vv}  \Big] } \dt
\\ &\geq \int_0^\tau \intQ{ \Big[ F(\mathbb{D}\vv) - \widehat{\mathbb{S}}:(\mathbb{D}\vv - \mathbb{D} \widehat{\vv}) - 
F(\mathbb{D} \widehat{\vv}) \Big] } \dt. 
\end{split}
\] 
Note that, strictly speaking as $\mathbb{D}\vv$ is only a measure we should interpret 
$$\int_0^\tau \!\!\int_{\Q} \widehat{\mathbb{S}} : \mathbb{D}\vv = \int_0^\tau \!\!\int_{\Q} \widehat{\mathbb{S}} : d\mathbb{D}\vv.$$
Moreover, as $\widehat{\mathbb{S}} \in \partial F(\mathbb{D} \widehat{\vv})$, we may infer that 
\[
\int_0^\tau \intQ{ \Big[ F(\mathbb{D}\vv) - \widehat{\mathbb{S}}:(\mathbb{D}\vv - \mathbb{D} \widehat{\vv}) - 
F(\mathbb{D} \widehat{\vv}) \Big] } \dt \geq 0.
\]
Indeed the inequality 
\[
\Big[ F(\mathbb{D}\vv) - \widehat{\mathbb{S}}:(\mathbb{D}\vv - \mathbb{D} \widehat{\vv}) - 
F(\mathbb{D} \widehat{\vv}) \Big] \geq 0
\]
can be derived by regularizing $\mathbb{D} \vv$ similarly to \cite[Lemma 3.2]{DemTema}.

Thus applying Gronwall's lemma to \eqref{R3}, we obtain the desired conclusion $\vv = \widehat \vv$, $\mathfrak{R}_v = 0$, and, finally, 
$\mathbb{S} = \widehat{\mathbb{S}}$. We have proved Theorem \ref{TM2}.

\section{Conditional regularity}
\label{A}

Our ultimate goal is to show Theorem \ref{TM3}. To this end observe that $\vv$ belonging to the regularity class \eqref{reg1} can be used 
as a test function in the momentum balance \eqref{D4}. After a routine manipulation, we obtain the energy balance in the form 
\begin{equation} \label{A1}
\intQ{ \frac{1}{2} |\vv|^2 (\tau, \cdot) } 
+ \int_0^\tau \intQ{ \mathbb{S} : \mathbb{D} \vv } \dt - \int_0^\tau \int_{\Q} \Grad \vv : \D \mathfrak{R}_v \dt  
= \intQ{ \frac{1}{2} |\vv_0|^2 }.
\end{equation}
Relation \eqref{A1} subtracted from the energy inequality \eqref{D5} gives rise to 
\begin{equation} \label{A2}
\int_{\Q} \D \frac{1}{2}{\rm trace}[\mathfrak{R}_v (\tau)]
+ \int_0^\tau \intQ{ \Big[ F(\mathbb{D} \vv) + F^*(\mathbb{S}) - \mathbb{S}: \mathbb{D} \vv \Big] } \dt  
\leq \int_0^\tau \int_{\Q} \Grad \vv : \D \mathfrak{R}_v \dt
\end{equation}
for a.a. $\tau \in (0,T)$. Consequently, by the standard Gronwall argument, we first deduce 
$\mathfrak{R}_v = 0$ and then 
\[
\mathbb{S} : \mathbb{D} \vv = F(\mathbb{D} \vv) + F^*(\mathbb{S}).
\]
We have proved Theorem \ref{TM3}.

\section{Concluding remarks}
\label{C}

The existence of local--in--time strong solutions has been shown by Bothe and Pr\"{u}ss \cite{BotPru} in the case of non--degenerate 
positive viscosity coefficient $\mu$. Similar result can be shown for the degenerate case $\mathbb{S} \equiv 0$ corresponding to the Euler system (note also that solutions of the Euler system are regular for regular initial data if $d = 2$). 
However, strictly speaking, the Euler system is excluded by hypothesis \eqref{DD5}. Obviously, the same approach works 
in this case as well, cf. Brenier, De Lellis, and Sz\' ekelyhidi \cite{BrDeSz}, Gwiazda et al. \cite{GSWW}, Sz\' ekelyhidi and Wiedemann \cite{SzeWie}. 
Although we could not find any relevant existence result in the degenerate mixed case (the fluid with activated viscosity according to Blechta, M\' alek, and Rajagopal 
\cite{BleMalRaj}), we believe the local--in--time existence of strong solution is in reach of the available analytic methods as soon as 
the viscosity coefficient $\mu$ is a sufficiently smooth function of $\mathbb{D}$.

For the sake of simplicity, we have also omitted the effect of external bulk force. It is easy to see that the latter can be accommodated in a straightforward manner.

Finally, we propose an extension of the method that accommodates both the inviscid (Euler) system as well as some viscous fluids with anisotropic viscous stress. Let 
\[
L: R^{d \times d}_{\rm sym} \to R^{d \times d}_{\rm sym}\ \mbox{be a linear mapping.} 
\]
For a convex function $F$ satisfying the hypotheses of Theorem \ref{TM1}, we consider 
\[
F_L (\mathbb{D}) = F( L \circ \mathbb{D} ),  
\]
which is again a convex function with ${\rm Dom} (F_L) = R^{d \times d}_{\rm sym}$. The associated conjugate function reads
\[
F^*_L = \mbox{(closure of)} \ \mathbb{S} \mapsto \inf_{\mathbb{M} \in R^{d \times d}_{\rm sym}} \left\{ F^*(\mathbb{M}),\ L^t \mathbb{M} = \mathbb{S} \right\}.
\]
The dissipative solutions are now defined exactly is in Definition \ref{DDD1}, where \eqref{D33} is replaced by 
\[
L [\mathbb{D} \vv] \in \mathcal{M}((0,T) \times \Q; R^{d \times d}_{\rm sym}), 
\]
and with $F_L$, $F^*_L$ in \eqref{D5}. Note that the choice $L \equiv 0$ gives rise to the Euler system, while 
$L = P$, where $P$ is a projection on a subspace of $R^{d \times d}_{\rm sym}$, corresponds to the anisotropic viscosity acting only in certain directions. The proofs of Theorems \ref{TM1}, \ref{TM2}, \ref{TM3} can be adapted in a direct manner.


\def\cprime{$'$} \def\ocirc#1{\ifmmode\setbox0=\hbox{$#1$}\dimen0=\ht0
  \advance\dimen0 by1pt\rlap{\hbox to\wd0{\hss\raise\dimen0
  \hbox{\hskip.2em$\scriptscriptstyle\circ$}\hss}}#1\else {\accent"17 #1}\fi}
\providecommand{\bysame}{\leavevmode\hbox to3em{\hrulefill}\thinspace}
\providecommand{\MR}{\relax\ifhmode\unskip\space\fi MR }
\providecommand{\MRhref}[2]{%
  \href{http://www.ams.org/mathscinet-getitem?mr=#1}{#2}
}
\providecommand{\href}[2]{#2}

\end{document}